\documentclass{article}
\usepackage{amsmath}
\usepackage[left=2.0cm,right=2.0cm,top=2.5cm,bottom=2.5cm]{geometry}

\numberwithin{equation}{section}
\newtheorem{thm}{Theorem}[section]
\newtheorem{cor}[thm]{Corollary}

\newtheorem{lem}[thm]{Lemma}

\newtheorem{ex}{Example}[section]

\newcommand{\be}{\begin{equation}}
\newcommand{\ee}{\end{equation}}
\newcommand{\ben}{\begin{enumerate}}
\newcommand{\een}{\end{enumerate}}

\newcommand{\qed}{\hspace*{\fill}Q.E.D.}  
\title{On Douglas general $(\alpha,\beta)$-metrics}
\author{Xiaoming Wang and Benling Li\footnote{Corresponding author. Research is supported by the NNSFC(11371209), ZPNSFC(LY13A010013) and K.C. Wong Magna Fund in Ningbo University.} \\
{\small \it Department of Mathematics, Ningbo University }\\
{\small \it Ningbo, Zhejiang Province 315211, P.R. China }\\
{\small \it wangxiaoming-0808@163.com;  libenling@nbu.edu.cn}
}

\begin{document}

\maketitle

\begin{abstract}
Douglas metrics are metrics with vanishing Douglas curvature which is an important projective invariant in Finsler geometry.
To find more Douglas metrics, in this paper we consider a class of Finsler
metrics called general $(\alpha,\beta)$-metrics, which are defined by a Riemannian metric $\alpha=\sqrt{a_{ij}(x)y^iy^j}$ and a $1$-form $\beta=b_i(x)y^i$. We obtain the differential equations that characterizes these metrics with vanishing Douglas curvature. By solving the equivalent PDEs, the metrics in this class are totally determined. Then many new Douglas metrics are constructed.


\end{abstract}

\section{Introduction}
Projective invariants play an important role in Finsler geometry. The  most used two projective invariants were introduced by J. Douglas in 1927
\cite{D1}. One of them is just Douglas curvature. A Finsler metric defined on an open subset in $R^n$ is called {\it Douglas metric} if its Douglas curvatures vanishes. Douglas curvature is a non-Riemannian quantity because it always vanishes for Riemannian metrics. In fact, all Riemannian metrics, Berwald metrics  and locally projectively flat Finsler metrics are special Douglas metrics. Thus, Douglas metrics form a rich class of metrics in Finsler geometry to show the difference and richness of Finsler geometry. Douglas metrics also have some applications in other problems.
Such as the application in constructing  non-Riemannian Einstein Finsler metrics from Douglas metrics \cite{C1,S1}.
Thus, it is natural to study and construct non-Rinmannian Douglas metrics.

  There are a lot of non-Riemannian metrics in Finsler geometry. Randers metric is the simplest non-Riemannian Finsler metric, which was introduced by the physicist G. Randers in \cite{R1}. The classification of Douglas Randers metric was obtained in \cite{B1}. As a generalization of Randers metric from the algebraic point of view, {\it $(\alpha,\beta)$-metrics} is defined by  the following form
$$F=\alpha \phi(s), \ \ \ s=\frac{\beta}{\alpha},$$
where $\alpha$ is a Riemannian metric, $\beta$ is a 1-form and $\phi (s)$ is a $C^{\infty}$ positive function. In \cite{L1}, we gave a characterization of Douglas $(\alpha,\beta)$-metrics with dimension $n \geq 3$.
A more general metric class called {\it general $(\alpha,\beta)$-metrics} was first introduced by C. Yu and H. Zhu in \cite{Y1}.
By definition, a {\it general $(\alpha,\beta)$-metric} is a Finsler metric expressed in the following form,
\be \label{gab} F=\alpha \phi(b^2,s), \ \ \ s=\frac{\beta}{\alpha},\ee
where $\alpha$ is a Riemannian metric, $\beta$ is a 1-form, $b:=\|\beta_x\|_\alpha$ and $\phi (b^2,s)$ is a $C^{\infty}$ positive function. It is easy to see that $(\alpha,\beta)$-metrics compose a special class in general $(\alpha,\beta)$-metrics. Another special class is defined by
 $\alpha$ being an Euclidean metric $|y|$ and $\beta$ being an inner product $\langle x,y\rangle$. In this case, the metric $F$ in (\ref{gab}) becomes a {\it spherically symmetric Finsler metric} in the following form
\be\label{E10}
F=|y| \phi(|x|^2,\frac{\langle x,y\rangle}{|y|}).
\ee
The related result about spherically symmetric Finsler metric can be found in \cite{MST}.

Recently, some non-Riemannian Einstein Finsler metrics were found in the class of Douglas $(\alpha,\beta)$-metrics \cite{C1,S1}. This motivates us to study the Douglas in the class of general $(\alpha,\beta)$-metrics. In \cite{LS1}, we classified projectively flat general $(\alpha,\beta)$-metrics when $\alpha$ is projectively flat.
In \cite{Zh1}, H. Zhu found a class of general $(\alpha,\beta)$-metrics with vanishing Douglas curvature under the condition that $\beta$ is closed and conformal with respect to $\alpha$, i.e., the covariant derivatives of $\beta$ with respect to $\alpha$ is $b_{i|j}=  c  a_{ij}$, where $c=c(x)\neq0$ is a scalar function on $M$.

In this paper, we remove the condition in \cite{Zh1} and find the equivalent equations of all the Douglas general $(\alpha,\beta)$-metrics. Based on these equivalent equations, we find some explicit new Douglas metrics.
We first give the following main theorem.

\begin{thm}\label{thm1.1} Let $F=\alpha \phi(b^2,\frac{\beta}{\alpha})$ be a non-Riemannian general $(\alpha,\beta)$-metric on an n-dimensional manifold with $n\geq3$. Suppose that $\beta$ is not parallel with respect to $\alpha$, then $F$ is a Douglas metric if and only if the function $\phi=\phi(b^2,s)$ satisfies the following PDE:
\be
\label{02}
\{b^3 [(1-c)s^2+cb^2]+[ ( \nu-\mu )s^2 - \nu b^2](b^2-s^2) \}\phi_{22}=2 b^5  (\phi_1-s \phi_{12})-[(\nu-\mu )s^2 - \nu b^2](\phi-s\phi_2)
\ee
and the covariant derivation of $\beta$ with respect to $\alpha$ satisfies the following equation:
\be b_{i|j}= kcb^2  a_{ij} + k(1-c) b_i b_j, \label{03}
\ee
where $k=k(x)$ is a scalar function,  $c=c(b^2)$,  $\mu=\mu(b^2)$ and $\nu=\nu(b^2)$
are  $C^{\infty}$ functions of $b^2$.
In this case,
$$G^i=\hat{G}^i + k\alpha \Big\{[(1-c)s^2+cb^2 ]\Theta +  b^2\Xi \Big\} y^i ,$$
where
\begin{eqnarray*}
\hat{G}^i & = & {^\alpha}G^i+\frac{k\alpha^2}{2b^3}\Big\{\nu b^2-(\nu-\mu) s^2\Big\} b^i,\\
\Theta &=&{(\phi-s \phi_2)\phi_2-s \phi \phi_{22} \over 2 \phi \Big[\phi-s \phi_2+(b^2-s^2)\phi_{22}\Big]},\\
\Xi & = &{b^2(\phi_1-s \phi_{12}) \phi_2 + s (b^2-s^2) \phi_1 \phi_{22} +s (\phi-s\phi_2)(2 \phi_1-s\phi_{12}) \over \phi \Big[\phi-s \phi_2+(b^2-s^2)\phi_{22}\Big]},
\end{eqnarray*}
here ${^\alpha}G^i$  are geodesic coefficients of  $\alpha$.
\end{thm}

Note that $\phi_1$ is the derivation of $\phi$ with respect to the first variable $b^2$ throughout the paper. The condition $n\geq3$ in the above theorem is natural. It is because that when the dimension $n=2$, obviously all Douglas metrics are just projectively flat metrics.
Obviously, when $c=1$, it is just the case in \cite{Zh1}. Since projectively flat Finsler metrics are special Douglas metrics, the PDE (\ref{02}) which is satisfied by $\phi$ is more general than the case in \cite{LS1}.

The above theorem tells us that there are many choices of the functions $k=k(x)$, $c=c(b^2)$, $\mu=\mu(b^2)$ and $\nu=\nu(b^2)$. To determine the Douglas metrics in Theorem \ref{thm1.1}, the efficient way is to solve the equivalent equations (\ref{02}) and (\ref{03}). Actually, we  obtain the following result about the general solutions of (\ref{02}).

\begin{thm}\label{thm1.2} Let $F=\alpha \phi(b^2,\frac{\beta}{\alpha})$ be a non-Riemannian general $(\alpha,\beta)$-metric on an n-dimensional manifold with $n\geq3$. If $\beta$ is not parallel with respect to $\alpha$, then the general solution of (\ref{02}) is given by
\be
\phi= s \Big\{h(b^2)- \xi(b^2) \int{\frac{ \Phi(\zeta(b^2,s))}{s^2 \sqrt{b^2-s^2}}ds}\Big\},
\label{04}
\ee
where
\be
\label{05}
\zeta(b^2,s)= \frac{b^2-s^2}{e^{ \int{ \frac{(1-c)b+\mu }{b^3}}db^2 }  + (b^2-s^2) \int{ \frac{\nu-\mu}{b^5} e^{ \int{ \frac{(1-c)b+\mu }{b^3}}db^2 }}db^2}
\ee
and
$$\xi(b^2)= e^{\int{\frac{1-c}{2 b^2}}db^2},$$
where $c=c(b^2)$, $h=h(b^2)$, $\mu=\mu(b^2)$ and $\nu=\nu(b^2)$ are $C^{\infty}$ functions of $b^2$.
$\Phi$ is a $C^{\infty}$ function of $\zeta$. To ensure the positive definite of $F$, $\Phi$ satisfies
\be
\label{lem511}
\frac{\Phi}{\sqrt{b^2-s^2}}>0, \ \ \ \ \ \Phi'\sqrt{b^2-s^2}>0
\ee
when $n\geq3$ or
\be
\label{lem512}
\Phi'\sqrt{b^2-s^2}>0
\ee
when $n=2$, where $\Phi'$ means the derivation of $\Phi$ with respect to the variable $\zeta$.
\end{thm}

Base on the above theorem, by choosing suitable functions $c(b^2)$,  $\mu(b^2)$, $\nu(b^2)$ and $\Phi(\zeta(b^2,s))$, we can construct some new Douglas metrics. Some examples are given in Section 6.

It is easy to see that spherically symmetric Finsler metric is a special class of general $(\alpha,\beta)$-metric. Let $\alpha=|y|$ and $\beta=\langle x,y\rangle$, then
 \[ b:=|x|, \ \ \ s:=\frac{\langle x,y\rangle}{|y|} \]
in (\ref{E10}).
Thus the following corollary is obvious by Theorem \ref{thm1.1} and Theorem \ref{thm1.2}.
\begin{cor}\label{cor1.1}
Let $F=|y| \phi(|x|^2,\frac{\langle x,y\rangle}{|y|})$ be  a non-Riemannian Finsler metric on an n-dimensional manifold with $n\geq3$, then $F$ is a Douglas metric if and only if $\phi$ satisfies
\be
\label{02cor}
[(\eta + f s^2)(b^2-s^2)-1 ]\phi_{22}+ 2(\phi_1-s \phi_{12}) + (\eta + f s^2)(\phi-s\phi_2)=0.
\ee
And  the general solution of (\ref{02cor}) is given by
$$\phi= s \Big\{\tilde{h}(b^2)-  \int{\frac{ \tilde{\Phi}(\tilde{\zeta}(b^2,s))}{s^2 \sqrt{b^2-s^2}}ds}\Big\},$$
where
$$\tilde{\zeta}(b^2,s)= \frac{b^2-s^2}{e^{ \int{ (f b^2 + \eta) }db^2 }  - (b^2-s^2) \int{ f e^{  \int{ (f b^2 + \eta) }db^2 }}db^2},$$
where
$f=f(b^2)$, $\eta=\eta(b^2)$ and $\tilde{h}=\tilde{h}(b^2)$   are arbitrary $C^{\infty}$ functions of $b^2$,
$\tilde{\Phi}$ is a $C^{\infty}$ function of $\tilde{\zeta}$. To ensure the positive definite of $F$, $\tilde{\Phi}$ satisfy
$$\frac{\tilde{\Phi}}{\sqrt{b^2-s^2}}>0, \ \ \ \ \ \tilde{\Phi}'\sqrt{b^2-s^2}>0$$
when $n\geq3$ or
$$
\tilde{\Phi}'\sqrt{b^2-s^2}>0
$$
when $n=2$, where $\tilde{\Phi}'$ means the derivation of $\tilde{\Phi}$ with respect to the variable $\tilde{\zeta}$.
 \end{cor}

The above corollary was first obtained  in \cite{MST} where  $\phi_1$  is the derivation of $\phi$ with respect to the  variable $b$. Essentially, the conclusion in \cite{MST} is the same as that of the Corollary \ref{cor1.1}.

\section{Preliminaries}
Let $M$ be a smooth $n$-dimensional manifold. A {\it Finsler metric} $F=F(x,y)$ on $M$ is a $C^\infty$ function, $F$: $TM$$\rightarrow$ $[0,\infty)$ with the following properties:
 (i) $F\geq0$ and $F(x,y)=0$ if and only if $y=0$; (ii)$F(x,\lambda y)=\lambda F(x,y)$ for all $\lambda>0$;
 (iii)$F$ is strongly convex, i.e., for any $y\neq0$, the matrix $g_{ij}=\frac{1}{2} [F^2]_{y^i y^j}$ is positive definite.
Specially, $F$ is called a {\it Riemannian metric} if $g_{ij}=g_{ij}(x)$. $F$ is  called a {\it Minkowskian metric} if $g_{ij}=g_{ij}(y)$.

The geodesics of a Finsler metric $F=F(x,y)$ on an open domain $\mathcal{U}\subset R^n$ can be defined by
\begin{equation*}
\frac{d^2 x^i}{dt^2}+2G^i(x,\frac{dx}{dt})=0,
\end{equation*}
where
\begin{equation*}
G^i = \frac{1}{4} g^{il}\Big{\{} [F^i]_{x^m y^l} y^m -[F^2]_{x^l} \Big{\}},
\end{equation*}
$g_{ij}:=\frac{1}{2}[F^2]_{y^i y^j}$ and $(g^{ij}):=(g_{ij})^{-1}$. The local functions $G^i=G^i(x,y)$ are called {\it geodesic coefficients}. \\

In \cite{D1}, J. Douglas introduced the quantity $D_y: T_xM \times T_xM \times T_xM \rightarrow T_xM $ which is a trillinear form $D_y(u,v,w)= D^i_{jkl}(y) u^jv^kw^l \frac{\partial}{\partial x^i}|_x$ defined by
$$D^i_{jkl}:=\frac{\partial^3}{\partial y^j \partial y^k \partial y^l} \Big( G^i-\frac{1}{n+1} \frac{\partial G^m}{\partial y^m} y^i \Big),$$
 $D^i_{jkl}$ are called the {\it Douglas curvature}. A Finsler metric $F$ is said to be a {\it Douglas metric} if $D^i_{jkl}=0$.
Base on this definition, by a direct computation,
Douglas metrics can  be characterized by
\be G^i=\frac{1}{2}\Gamma^i_{jk}(x)y^j y^k+P(x,y)y^i,\label{11}
\ee
where $\Gamma^i_{jk}(x)$ are local functions on $M$ and $P(x,y)$ is a local positively homogeneous function of degree one  \cite{Sh2}.\\
Specially, $F$ is called {\it locally projectively flat} if $G^i = P y^i$. 

In this paper, we consider the general $(\alpha,\beta)$-metrics. The following lemma was proved in  \cite{Y1}.
\begin{lem}\label{lem2.1} {\rm (\cite{Y1})} Let $M$ be an n-dimensional manifold. $F=\alpha \phi(b^2,\frac{\beta}{\alpha})$ is a Finsler metric on $M$ for any Riemannian metric $\alpha$ and 1-form $\beta$ with $\|\beta\|_\alpha < b_0$ if and only if $\phi=\phi(b^2,s)$ is a positive $C^{\infty}$ function satisfying
\be \phi-s \phi_2 >0, \ \  \phi-s \phi_2 + (b^2-s^2)\phi_{22}>0,\label{13}
\ee
when $n\geq3$ or $$\phi-s \phi_2 + (b^2-s^2)\phi_{22}>0,$$
when $n=2$, where $s$ and $b$ are arbitrary numbers with $|s|\leq b <b_0.$
\end{lem}

Denote the covariant derivative of the 1-form $\beta$ with respect to the  Riemannian metric $\alpha$ by $b_{i|j}$. Moreover, for simplicity, let
$$r_{ij}=\frac{1}{2}(b_{i|j}+b_{j|i}), \ s_{ij}=\frac{1}{2}(b_{i|j}-b_{j|i}), \  r_{00}= r_{ij}y^iy^j, \ s^i_{ \ 0}=a^{ij}s_{jk}y^k,$$
$$r_i= b^j r_{ji}, \ s_i=b^js_{ji}, \ r_0 = r_iy^i, \ s_0 =s_iy^i, \ r^i=a^{ij}r_j, \ s^i=a^{ij}s_j, \ r=b^ir_i.$$

The geodesic coefficients $G^i$ of a general $(\alpha,\beta)$-metric $F=\alpha \phi(b^2,\frac{\beta}{\alpha})$ were given in \cite{Y1} as the following
\begin{equation}\label{14}
\begin{split}&G^i={^\alpha}G^i+{\alpha}Q s^i_{\ 0}+\Big\{\Theta(-2 \alpha Q s_0+r_{00}+2 \alpha^2 R r)+\alpha\Omega(r_0+s_0)\Big\}\frac{y^i}{\alpha}\\&
\ \ \ \ \ \ +\Big\{\Psi(-2\alpha Q s_0+r_{00}+2\alpha^2 R r)+\alpha \Pi(r_0+s_0)\Big\}b^i-\alpha^2 R (r^i+s^i),
\end{split}
\end{equation}
where
$$Q=\frac{\phi_2}{\phi-s \phi_2},\ \ \Theta=\frac{(\phi-s \phi_2)\phi_2-s \phi \phi_{22}}{2 \phi [\phi-s \phi_2+(b^2-s^2)\phi_{22}]},\ \ \Psi=\frac{\phi_{22}}{2 [\phi-s \phi_2+(b^2-s^2)\phi_{22}]},$$
$$R=\frac{\phi_1}{\phi-s \phi_2}, \ \ \Pi=\frac{(\phi-s \phi_2)\phi_{12}-s \phi_1 \phi_{22}}{(\phi-s \phi_2) [\phi-s \phi_2+(b^2-s^2)\phi_{22}]},\ \ \Omega=\frac{2 \phi_1}{\phi}-\frac{s \phi+(b^2-s^2)\phi_2}{\phi}\Pi.$$
Here ${^\alpha}G^i$ are geodesic coefficients of $\alpha$.

By (\ref{11}),\ Douglas metrics can be also characterized by the following equations \cite{B1}:
\be
G^i y^j-G^j y^i=\frac{1}{2}(\Gamma^i_{kl} y^j-\Gamma^j_{kl} y^i)y^k y^l.
\label{15}
\ee
Substituting (\ref{14}) into (\ref{15}), it  can be easily shown that a general $(\alpha,\beta)$-metric is a Douglas metric if and only if
 \begin{equation}\label{16}
\begin{split}&\{\Psi(-2\alpha Q s_0+r_{00}+2\alpha^2 R r)+\alpha\Pi(r_0+s_0)\}(b^i y^j-b^j y^i)
-\alpha^2 R\{(r^i+s^i)y^j-(r^j+s^j)y^i\}\\&+\alpha Q(s^i_{\ 0} y^j-s^j_{\ 0}y^i)
=\frac{1}{2}(G^i_{kl}y^j-G^j_{kl}y^i)y^k y^l,
\end{split}
\end{equation}
where $G^i_{kl}:=\Gamma^i_{kl}-\gamma^i_{kl}$ and $\gamma^i_{kl}:=\frac{\partial^2 G^i_\alpha}{\partial y^k\partial y^l}$. \\

The following lemmas are needed in Section 4 and Section 5.
\begin{lem}\label{lem2.2} If $Q=\iota_1 s$, where $\iota_1=\iota_1(b^2)$ is independent of $s$, then
\be
\label{lem2.21}
\phi(b^2,s)=\iota_2 \sqrt{1+\iota_1 s^2},
\ee
where $\iota_2=\iota_2(b^2)$ is a function of $b^2$.
\end{lem}
{\bf Proof}: By a direct computation, we have (\ref{lem2.21}).

\qed
\begin{lem}\label{lem2.3} If $\Psi=\iota_3+ \frac{\iota_4 s^2}{b^2-s^2}$, where $\iota_3=\iota_3(b^2)$ and $\iota_4=\iota_4(b^2)$ are independent of $s$, then
\be
\label{lem2.31}
\phi(b^2,s)=s\Big\{\iota_6-\int{ \frac{(b^2-s^2)^{-\frac{b^2\iota_4}{2b^2\iota_4-1}}[2(\iota_4-\iota_3)s^2+2\iota_3b^2-1]^{\frac{1}{2(2b^2\iota_4-1)}}\iota_5}{s^2}ds} \Big\},
\ee
where $\iota_5=\iota_5(b^2)$ and $\iota_6=\iota_6(b^2)$  are functions of $b^2$.
\end{lem}
{\bf Proof}: By the assumption, we have
\be
\label{lem2.32}
(b^2-s^2)[1-2(\iota_4-\iota_3)s^2-2\iota_3b^2]\phi_{22}=2[(\iota_4-\iota_3)s^2+\iota_3b^2](\phi-s\phi_2).
\ee
Note that $(\phi-s \phi_2)_2 = - s \phi_{22}$.
By setting
\be
\label{lem2.34}\varphi:=\phi-s\phi_2,
\ee
  (\ref{lem2.32}) is equivalent to
\be
\label{lem2.33}
(b^2-s^2)[2(\iota_4-\iota_3)s^2+2\iota_3b^2-1] \varphi_2=2s[(\iota_4-\iota_3)s^2+\iota_3b^2] \varphi.
\ee
By a direct computation,  the solution of (\ref{lem2.33}) is
$$\varphi=(b^2-s^2)^{-\frac{b^2\iota_4}{2b^2\iota_4-1}}[2(\iota_4-\iota_3)s^2+2\iota_3b^2-1]^{\frac{1}{2(2b^2\iota_4-1)}}\iota_5,$$
where $\iota_5=\iota_5(b^2)$ is a function of $b^2$. Plugging the above equation into (\ref{lem2.34}), we have
\be
\label{lem2.35}
\phi-s\phi_2=(b^2-s^2)^{-\frac{b^2\iota_4}{2b^2\iota_4-1}}[2(\iota_4-\iota_3)s^2+2\iota_3b^2-1]^{\frac{1}{2(2b^2\iota_4-1)}}\iota_5.
\ee
Let
\be
\phi=s \theta,
\label{lem2.36}
\ee
then
\be
\phi-s \phi_2 = - s^2 \theta_2.
\label{lem2.37}
\ee
By (\ref{lem2.35}) and (\ref{lem2.37}), we obtain
\be
\theta= \iota_6-\int{ \frac{(b^2-s^2)^{-\frac{b^2\iota_4}{2b^2\iota_4-1}}[2(\iota_4-\iota_3)s^2+2\iota_3b^2-1]^{\frac{1}{2(2b^2\iota_4-1)}}\iota_5}{s^2}ds}
\label{lem2.38}
\ee
for some $C^{\infty}$ functions $\iota_6=\iota_6(b^2)$.
Hence, plugging (\ref{lem2.38}) into (\ref{lem2.36}) gives (\ref{lem2.31}).

\qed

\section{Sufficient Conditions of Theorem \ref{thm1.1}}
In this section, we are going to prove that the sufficient conditions for a general $(\alpha,\beta)$-metric to be a Douglas metric.
\begin{lem}\label{lem3.1}
Let $F=\alpha \phi(b^2,\frac{\beta}{\alpha})$ be a non-Riemannian general $(\alpha,\beta)$-metric on an n-dimensional manifold with $n\geq2$. Suppose that $\beta$ satisfies
\be b_{i|j}= kcb^2  a_{ij} + k(1-c) b_i b_j \label{lem03case1}
\ee
and $\phi=\phi(b^2,s)$ satisfies the following PDE:
\be
\label{lem02case1}
\{b^3 [(1-c)s^2+cb^2]+[ ( \nu-\mu )s^2 - \nu b^2](b^2-s^2) \}\phi_{22}=2 b^5  (\phi_1-s \phi_{12})-[(\nu-\mu )s^2 - \nu b^2](\phi-s\phi_2),
\ee
where $k=k(x)$ is a scalar function,  $c=c(b^2)$,  $\mu=\mu(b^2)$ and $\nu=\nu(b^2)$
are arbitrary $C^{\infty}$ functions of $b^2$.
Then the geodesic coefficients $G^i=G^i(x,y)$ of $F$ are given by
\be
\label{lem16}
G^i=\hat{G}^i + k\alpha \Big\{[(1-c)s^2+cb^2]\Theta + b^2\Xi \Big\} y^i ,
\ee
where
\begin{eqnarray*}
\bigskip
\hat{G}^i & = & {^\alpha}G^i+\frac{k\alpha^2}{2b^3}\Big\{\nu b^2-(\nu-\mu) s^2\Big\} b^i\\
\Theta &=&{(\phi-s \phi_2)\phi_2-s \phi \phi_{22} \over 2 \phi \Big[\phi-s \phi_2+(b^2-s^2)\phi_{22}\Big]}\\
\Xi & = &{b^2(\phi_1-s \phi_{12}) \phi_2 + s (b^2-s^2) \phi_1 \phi_{22} +s (\phi-s\phi_2)(2 \phi_1-s\phi_{12}) \over \phi \Big[\phi-s \phi_2+(b^2-s^2)\phi_{22}\Big]}.
\bigskip
\end{eqnarray*}
Here ${^\alpha}G^i$ are geodesic coefficients of $\alpha$.
\end{lem}
{\bf Proof}: By the assumption (\ref{lem03case1}), we have
$$ r_{00}= kcb^2  \alpha^2 + k(1-c) \beta^2, \ \ \ r_0 = kb^2 \beta, \ \ \ r^i= kb^2 b^i, \ \ \ r=k b^4, \ \ \ s_0 = 0, \ \ \ s^i_{ \ 0} = 0.$$
Substituting them into (\ref{14}) yields
\begin{equation}\label{lem14}
G^i={^\alpha}G^i+k\alpha \{\Theta[(1-c)s^2+cb^2+2 b^4 R ]+b^2 s\Omega \}y^i +k \alpha^2 \{\Psi[(1-c)s^2+cb^2+2 b^4 R ]+b^2s\Pi-b^2 R \}b^i.
\end{equation}
By substituting the expression of $\Theta, R, \Psi, \Pi$ and $\Omega$ into the above equation we have
\begin{equation}\label{lem15}
\begin{split}&G^i={^\alpha}G^i+k\alpha \Big\{[(1-c)s^2+cb^2 ]{(\phi-s \phi_2)\phi_2-s \phi \phi_{22} \over 2 \phi [\phi-s \phi_2+(b^2-s^2)\phi_{22}]} \\&
\ \ \ \ \ \ \ + b^2{b^2(\phi_1-s \phi_{12}) \phi_2 + s (b^2-s^2) \phi_1 \phi_{22} +s (\phi-s\phi_2)(2 \phi_1-s\phi_{12}) \over \phi [\phi-s \phi_2+(b^2-s^2)\phi_{22}]} \Big\}y^i\\&
\ \ \ \ \ \ \ + k\alpha^2 {[(1-c)s^2+cb^2 ]\phi_{22}- 2 b^2 (\phi_1-s \phi_{12}) \over 2 [\phi-s\phi_2+(b^2-s^2)\phi_{22}]}b^i.
\end{split}
\end{equation}

On the other hand, by (\ref{lem02case1}), we obtain
\be
\label{lem17}
\frac{[(1-c)s^2+cb^2 ]\phi_{22}- 2 b^2 (\phi_1-s \phi_{12}) }  {\phi-s\phi_2+(b^2-s^2)\phi_{22}}=\frac{\nu b^2-(\nu-\mu) s^2}{ b^3}.
\ee
Plugging (\ref{lem17}) into (\ref{lem15}) yields (\ref{lem16}).

\qed\\

{\bf Proof of the sufficiency of Theorem \ref{thm1.1}:}
Note that $\hat{G}^i$  are quadratic in $y\in T_xM$. Thus $F$ is a Douglas metric by Lemma \ref{lem3.1} and (\ref{11}).

\qed

The proof of the converse is given in the following two sections. Firstly, we prove $\beta$ is closed in Section 4. Then the PDE of $\phi$ can be obtained in Section 5.

\section{$\beta$ is closed}
In this section, we mainly prove that the 1-form $\beta$ is closed for Douglas metrics.
In order to analysis (\ref{16}), we choose a special coordinate system at a point as in \cite{Sh1}. Take an orthonormal basis at any fixed point $x_0$  such that
\be\alpha=\sqrt{\sum(y^i)^2}, \ \ \ \ \  \beta=b y^1,\label{20}
\ee
where $b:=\|\beta_x\|_\alpha$.
Make a change of coordinates: $(s,y^a)\rightarrow(y^i)$ by
$$y^1=\frac{s}{\sqrt{b^2-s^2}}\bar{\alpha},\ \ \ \ \ y^a=y^a,\ \ a = 2, \ldots, n $$
where $\bar{\alpha}:=\sqrt{\sum^n_{a=2}(y^a)^2}$. Then
\be\alpha=\frac{b}{\sqrt{b^2-s^2}}\bar{\alpha},\ \ \ \ \beta=\frac{b s}{\sqrt{b^2-s^2}}\bar{\alpha}.
\label{30}
\ee
Let
$$\bar{r}_{10}:=\sum^n_{a=2} r_{1a} y^a,\ \  \bar{s}_{10}:=\sum^n_{a=2} s_{1a} y^a,\ \ \bar{r}_{00}:=\sum^n_{a,b=2} r_{ab} y^a y^b,\ \ \bar{s}^1_{\ 0}:=\sum^n_{a=2} s^1_{\ a} y^a,\ \ \bar{s}_{0}:=\sum^n_{a=2} s_{a} y^a,\ \ \bar{r}_{0}:=\sum^n_{a=2} r_{a} y^a,$$
and so on. It's easy to see that $s_1 = b s^1_{\ 1} = 0,\ r_1 = b r_{11}, \ \bar{s}_{0} = b \bar{s}_{10}, \ \bar{r}_{0} = b \bar{r}_{10}.$
Thus we have
$$s^1_{\ 0}=\bar{s}^1_{\ 0},\ \ \ \ \ s^a_{\ 0}=\frac{s \bar{s}^a_{\ 1} }{\sqrt{b^2-s^2}} \bar{\alpha}+\bar{s}^a_{\ 0}, \ \ \ \  \ s_0 = \bar{s}_0,\ \ \ \ \ r_0=\frac{b s r_{11}}{\sqrt{b^2-s^2}} \bar{\alpha}+b \bar{r}_{10},$$
$$r_{00}=r_{11} \frac{s^2}{b^2-s^2} \bar{\alpha}^2+2 \bar{r}_{10} \frac{s}{\sqrt{b^2-s^2}} \bar{\alpha}+\bar{r}_{00},\ \ \ \ \ \ r=b^2 r_{11}.$$
Let
$$\bar{G}^a_{10}=G^a_{1b} y^b,\ \ \ \ \ \bar{G}^a_{01}=G^a_{b1} y^b,\ \ \ \ \ \bar{G}^a_{00}=G^a_{bc} y^b y^c.$$
Then
$$G^a_{kl}y^k y^l=G^a_{11} \frac{s^2}{b^2-s^2} \bar{\alpha}^2+ \bar{G}^a_{10} \frac{s}{\sqrt{b^2-s^2}} \bar{\alpha}+\bar{G}^a_{01} \frac{s}{\sqrt{b^2-s^2}} \bar{\alpha}+\bar{G}^a_{00}.$$
Plugging the above identities into (\ref{16}), for\ $i=1$, $j=a$,  we get a system of equations in the following form
$$A\bar{\alpha}^3+B\bar{\alpha}^2+C\bar{\alpha}+E=0,$$
where $A$, $B$, $C$ and $E$ are polynomials in $y^a$.
By $\bar{\alpha}=\sqrt{\sum^n_{a=2}(y^a)^2}$ is a irrational function of $y$, we have
$$A\bar{\alpha}^2+C=0, \ \ \ \ B\bar{\alpha}^2+E=0.$$
They are equivalent to
\begin{equation}\label{21}
\begin{split}&\frac{\bar{\alpha}^2}{b^2-s^2}\{b Q \bar{s}^a_{\ 0}s-b r_{11} (\Psi s^2+2 \Psi R b^4 + \Pi b^2 s - R b^2) y^a \}-\Psi \bar{r}_{00}b y^a\\&=\frac{1}{2(b^2-s^2)}\{(\bar{G}^a_{10}+\bar{G}^a_{01})s^2-G^1_{11}s^2 y^a\}\bar{\alpha}^2-\frac{1}{2}\bar{G}^1_{00}y^a
\end{split}
\end{equation}
and
\begin{equation}\label{22}
\begin{split}&\frac{\bar{\alpha}^2}{b^2-s^2}\{b Q s^a_{\ 1} s^2-b^2 s R(r_a+s_a)\}+\{2\Psi Q b^2 \bar{s}_{10}-2\Psi s \bar{r}_{10}-\Pi b^2(\bar{r}_{10}+\bar{s}_{10})-Q s^1_{\ 0}\}b y^a\\&=\frac{G^a_{11} s^3}{2(b^2-s^2)}\bar{\alpha}^2+\frac{1}{2}\{\bar{G}^a_{00}-(\bar{G}^1_{10}+\bar{G}^1_{01})y^a\}s .
\end{split}
\end{equation}
Similarly, for\ $i=a$, $j=b$ $(a\neq b)$,  we get
\begin{equation}\label{23}
\begin{split}&\frac{b s}{b^2-s^2}(s^a_{\ 1}y^b-s^b_{\ 1}y^a)Q \bar{\alpha}^2+\frac{b^2 R}{b^2-s^2}\{(r_b+s_b)y^a-(r_a+s_a)y^b\}\bar{\alpha}^2\\&=
\frac{s^2}{2(b^2-s^2)}(G^a_{11}y^b-G^b_{11}y^a)\bar{\alpha}^2+\frac{1}{2}(\bar{G}^a_{00}y^b-\bar{G}^b_{00}y^a),
\end{split}
\end{equation}

\be (\bar{s}^a_{\ 0}y^b-\bar{s}^b_{\ 0}y^a)b Q=\frac{s}{2}\{(\bar{G}^a_{10}+\bar{G}^a_{01})y^b-(\bar{G}^b_{10}+\bar{G}^b_{01})y^a\}. \label{24}
\ee

Now we can prove the following
\begin{lem}\label{lem4.1} If a non-Riemannian general $(\alpha,\beta)$-metric $F=\alpha \phi(b^2,\frac{\beta}{\alpha})$ is a Douglas metric on an n-dimensional manifold with $n\geq 3$, then $\beta$ is closed.
\end{lem}
{\bf Proof}: Letting $s=0$ in (\ref{23}), we get
$$\frac{1}{2}(\bar{G}^a_{00}y^b-\bar{G}^b_{00}y^a)=R(b^2,0) \{(r_b+s_b)y^a-(r_a+s_a)y^b\}\bar{\alpha}^2,$$
Substituting it back into (\ref{23}) yields
\be b s (s^a_{\ 1}y^b-s^b_{\ 1}y^a)Q+ [b^2 R-R(b^2,0)(b^2-s^2)] \{(r_b+s_b)y^a-(r_a+s_a)y^b\}= \frac{s^2}{2}(G^a_{11}y^b-G^b_{11}y^a). \label{25}
\ee
By the arbitrary of $y^a$ and $y^b$, there exist $y^a_0$ and $y^b_0$ such that
$$(r_b+s_b)y^a_0 =(r_a+s_a)y^b_0.$$
Then (\ref{25}) becomes
\be b (s^a_{\ 1}y^b_0-s^b_{\ 1}y^a_0)Q= \frac{s}{2}(G^a_{11}y^b_0-G^b_{11}y^a_0). \label{26}
\ee
If $Q=\iota_1 s$, where $\iota_1=\iota_1(b^2)$ is independent of $s$, by Lemma \ref{lem2.2}, we have  $F=\iota_2 \sqrt{\alpha^2+\iota_1 \beta^2}$ is of Riemannian type, where $\iota_2=\iota_2(b^2)$ is a function of $b^2$.
This case is exclude in our assumption of Lemma \ref{lem4.1}.
Thus we conclude that $Q \neq \iota_1 s$. Since $b\neq0$, by (\ref{26}) we obtain
$$s^a_{\ 1}y^b_0-s^b_{\ 1}y^a_0=G^a_{11}y^b_0-G^b_{11}y^a_0=0.$$
Then by the arbitrary choice of $x_0$, we obtain
\be \frac{s_a}{s_b}=\frac{b s^a_{\ 1}}{b s^b_{\ 1}}=\frac{r_a+s_a}{r_b+s_b}=\frac{r_a}{r_b}=\frac{G^a_{11}}{G^b_{11}}. \label{27}
\ee
In the above equation, if the denominator is  zero, which implies that numerator is  zero too.
Set
\be
\label{c1}
s^a_{\ 1}=c_1(r_a+s_a)
\ee
and
\be
\label{c2}
G^a_{11}=c_2(r_a+s_a),
\ee
where $c_1=c_1(x)$ and $c_2=c_2(x)$  are two numbers at the point $x_0$. Substituting (\ref{c1}) and (\ref{c2}) into (\ref{25}) yields
\be
\label{200}
\{(r_a+s_a)y^b-(r_b+s_b)y^a\}\{ b s c_1 Q-b^2 R -\frac{c_2}{2}s^2 + c_3(b^2-s^2)\}=0,
\ee
where $c_3=R(b^2,0)$ is a function of $b^2$.

We claim that $(r_a+s_a)y^b-(r_b+s_b)y^a = 0$.
If $(r_a+s_a)y^b-(r_b+s_b)y^a \neq 0$,
then by (\ref{200}), we have
\be
\label{28}
 b s c_1 Q-b^2 R=\frac{c_2}{2}s^2 - c_3(b^2-s^2).
\ee
Differentiating (\ref{28}) with respect to $s$ for three times yields
$$c_1 b (s Q)_{sss}-b^2 (R)_{sss}=0.$$
Then $c_1$ is a function of $b^2$. Thus by (\ref{28}), $c_2$  is also a function of $b^2$.
Substituting the expressions of $Q$ and $R$ into (\ref{28}), because $\phi-s \phi_2 > 0$, we get
\be 2b s c_1 \phi_2- 2b^2 \phi_1= [(c_2+2c_3) s^2-2b^2c_3](\phi-s \phi_2). \label{29}
\ee
By defining
\be
\phi =\omega s e^{\int{\frac{c_1}{b}}db^2},
\label{2}
\ee
we can prove $\omega$ is an even function in $s$.
Substituting (\ref{2}) into (\ref{29}), which can be rewritten  as
\be
-2 b^2 \omega_1 +  [(2 b c_1-2b^2 c_3) s + (c_2+2c_3) s^3]\omega_2 = 0.
\label{3}
\ee
We are going to get the general solution of (\ref{3}).
The characteristic equation of  (\ref{3}) is
\be \frac{db^2}{-2 b^2} = \frac{d s}{(2 b c_1-2b^2 c_3) s + (c_2+2c_3) s^3},
\label{4}
\ee
which is equivalent to
\be  \frac{d s}{db^2} = -\frac{c_1-bc_3}{b} s - \frac{c_2+2c_3}{2 b^2}s^3.
\label{5}
\ee
Obviously, it is a Bernoulli equation. Consider the following variable substitution
\be \chi(b^2) = \frac{1}{s^2(b^2)}.
\label{6}
\ee
Then (\ref{5}) can be rewritten as
\be \frac{d \chi}{db^2} = \frac{2 c_1-2 b c_3}{b} \chi + \frac{c_2+2c_3}{ b^2}.
\label{8}
\ee
This is a linear 1-order ODE of $\chi$, whose solution is
\be \chi =e^{2 \int{ \frac{c_1- b c_3}{b}}db^2 }\Big[\lambda+\int{ \frac{c_2+2c_3}{b^2} e^{ -2\int{ \frac{c_1- b c_3}{b}}db^2 }}db^2 \Big],
\label{9}
\ee
where $\lambda$ is an arbitrary constant.
Then by (\ref{6}) and (\ref{9}) we get
\be
\label{7}
s=s(b^2)=\sqrt{\frac{e^{-2 \int{ \frac{c_1- b c_3}{b}}db^2 }}{\lambda +\int{ \frac{c_2+2c_3}{b^2} e^{ -2\int{ \frac{c_1- b c_3}{b}}db^2 }}db^2}}.
\ee
$s=0$ is also a solution of (\ref{5}), it is excluded.
It follows from (\ref{7}) that
$$\frac{s^2}{  e^{ -2\int{ \frac{c_1-bc_3}{b}}db^2 }-s^2 \int{ \frac{c_2+2c_3}{b^2} e^{ -2\int{ \frac{ c_1-bc_3}{b}}db^2 } } db^2}=  \frac{1}{\lambda}.$$
Hence the solution of (\ref{3}) is
\be
\omega = \Upsilon \Big(\frac{s^2}{  e^{ -2\int{ \frac{c_1-bc_3}{b}}db^2 }-s^2 \int{ \frac{c_2+2c_3}{b^2} e^{ -2\int{ \frac{ c_1-bc_3}{b}}db^2 } } db^2} \Big),
\label{10}
\ee
where $\Upsilon$ is any differentiable function. Plugging (\ref{10}) into (\ref{2}) gives
\be
\phi=s e^{\int{\frac{c_1}{b}}db^2}\Upsilon \Big(\frac{s^2}{  e^{ -2\int{ \frac{c_1-bc_3}{b}}db^2 }-s^2 \int{ \frac{c_2+2c_3}{b^2} e^{ -2\int{ \frac{ c_1-bc_3}{b}}db^2 } } db^2} \Big).
\label{111}
\ee
Obviously, $\Upsilon$ is an even function in $s$ and $s$ is an odd function, then $\phi$ is an odd function in $s$.
Because $\phi$ is also a $C^{\infty}$ positive function, it must be meaningful at the origin, which means that $\phi=0$ when $s=0$. It is contradict to $\phi(b^2,0) \neq 0$. Because when $\beta=0$, $F=\alpha\phi(b^2,0)$ should be a Riemannian metric. Thus  $$(r_a+s_a)y^b-(r_b+s_b)y^a = 0.$$
By the arbitrary of $y^a$ and $y^b$, we have
\be r_a+s_a=0.
\label{213}
\ee
Therefore, by  (\ref{27}) and (\ref{213}), we obtain
\be
\label{214}
s^a_{\ 1}=0,
\ee
\be
\label{215}
 r_{1a}=0
 \ee
and
$$ G^a_{11}=0.$$

Differentiating (\ref{24}) with respect to $y^c$ ($c\neq a$ ,$c\neq b$) and $y^a$ ($a\neq b$) yields
\be
\label{217}
b s^b_{\ c}  Q = \frac{s}{2} (G^b_{1c}+G^b_{c1}).
\ee
Since $Q \neq \iota_1 s$, where $\iota_1=\iota_1(b^2)$ is independent of $s$, we conclude that
\be s^a_{\ b}=0. \label{216}
\ee
Plugging (\ref{216}) back  into (\ref{217}) yields
\be
\label{218}
G^b_{1c}+G^b_{c1}=0, \ \ (b\neq c).
\ee
Substituting (\ref{216}) and (\ref{218}) into (\ref{24}), we also get
$$ \bar{G}^a_{10}+\bar{G}^a_{01}=Hy^a,$$
where $H=H(x)=G^a_{1a}+G^a_{a1} $ is a number at the point $x_0$.

Thus by (\ref{214}) and (\ref{216}), we obtain
$s_{ij}=0$ which means that $\beta$ is closed.

\qed

\section{Necessary Conditions of Theorem \ref{thm1.1}}
In this section, we are going to prove the necessity of Theorem \ref{thm1.1}. It can be obtained by the following lemma.
\begin{lem}\label{lem5.1} Let $F=\alpha \phi(b^2,\frac{\beta}{\alpha})$ be a non-Riemannian general $(\alpha,\beta)$-metric on an n-dimensional manifold with $n\geq3$. Assume that $\beta$ is not parallel with respect to $\alpha$. If $F$ is a Douglas metric, then  there are four scalar functions $k = k (x) $,  $c = c(b^2)$, $\mu=\mu(b^2)$ and $\nu=\nu(b^2)$ such that
\be
\label{lem4.102case1}
\{b^3 [(1-c)s^2+cb^2]+[ ( \nu-\mu )s^2 - \nu b^2](b^2-s^2) \}\phi_{22}=2 b^5  (\phi_1-s \phi_{12})-[(\nu-\mu )s^2 - \nu b^2](\phi-s\phi_2)
\ee
and
\be b_{i|j}= kcb^2  a_{ij} + k(1-c) b_i b_j . \label{lem4.103case1}
\ee
\end{lem}
{\bf Proof}:
By Lemma \ref{lem4.1}, we get
$$s_{ij}=0,\ \ \ r_{1a}=0, \ \ \ G^a_{11}=0,\ \ \ \bar{G}^a_{10}+\bar{G}^a_{01}=Hy^a.$$
Then (\ref{21})  is reduced to
\begin{equation}\label{31}
2\bar{\alpha}^2 b r_{11} (\Psi s^2+2 \Psi R b^4 + \Pi b^2 s - R b^2)+2 \Psi \bar{r}_{00}b(b^2-s^2)=(G^1_{11}-H)\bar{\alpha}^2 s^2+\bar{G}^1_{00}(b^2-s^2).
\end{equation}
Letting $s=0$ in (\ref{31}), we have
\be  \bar{G}^1_{00} = 2\bar{\alpha}^2 b r_{11} [2 \Psi (b^2,0) R(b^2,0) b^2  - R(b^2,0)] + 2 \Psi(b^2,0) \bar{r}_{00}b.
\label{32}
\ee
Plugging (\ref{32}) into (\ref{31}), we get
\begin{equation}\label{322}
\begin{split}&\Big\{2 b r_{11} \{\Psi s^2+2 \Psi R b^4 + \Pi b^2 s - R b^2-[2 \Psi (b^2,0) R(b^2,0) b^2  - R(b^2,0)](b^2-s^2)\}+(H-G^1_{11})s^2\Big\}\bar{\alpha}^2\\&=
2 b(b^2-s^2)\{ \Psi(b^2,0)- \Psi\}\bar{r}_{00}.
\end{split}
\end{equation}
Differentiating (\ref{322}) with respect to $y^a$ and $y^b$ yields
\begin{equation}\label{323}
\begin{split}&\Big\{2 b r_{11} \{\Psi s^2+2 \Psi R b^4 + \Pi b^2 s - R b^2-[2 \Psi (b^2,0) R(b^2,0) b^2  - R(b^2,0)](b^2-s^2)\}+(H-G^1_{11})s^2\Big\}\delta_{ab}\\&=
2 b(b^2-s^2)\{ \Psi(b^2,0)- \Psi\}r_{ab}.
\end{split}
\end{equation}
Then  by (\ref{323}), there exists a scalar function $\lambda=\lambda(x)$ such that
\be r_{ab}= \lambda \delta_{ab}. \label{33}
\ee
We may set
\be
r_{11} = k b^2, \ \ \ (k=k(x)).
\label{34case1}
\ee
If $r_{11}=0$, plugging it into (\ref{323}), we have $\Psi=\iota_3+ \frac{\iota_4 s^2}{b^2-s^2}$, where $\iota_3=\iota_3(b^2)$ and $\iota_4=\iota_4(b^2)$ are independent of $s$. By Lemma \ref{lem2.3}, $\phi$ is given by
$$\phi(b^2,s)=s\Big\{\iota_6-\int{ \frac{(b^2-s^2)^{-\frac{b^2\iota_4}{2b^2\iota_4-1}}[2(\iota_4-\iota_3)s^2+2\iota_3b^2-1]^{\frac{1}{2(2b^2\iota_4-1)}}\iota_5}{s^2}ds} \Big\},$$
where $\iota_5=\iota_5(b^2)$ and $\iota_6=\iota_6(b^2)$  are functions of $b^2$.
Then
$$\phi_s=\iota_6-\int{\frac{\varphi}{s^2} ds }- \frac{\varphi}{s},$$
where $$\varphi=(b^2-s^2)^{-\frac{b^2\iota_4}{2b^2\iota_4-1}}[2(\iota_4-\iota_3)s^2+2\iota_3b^2-1]^{\frac{1}{2(2b^2\iota_4-1)}}\iota_5.$$
It is easy to see that $\phi_s$ is not a $C^{\infty}$ function. Then it is excluded. Thus $r_{11} \neq 0$, which means $k\neq0$. Then there exists a scalar function $c=c(x)$ such that
\be
\label{324}
\lambda=ckb^2.
\ee
Plugging (\ref{324}) into (\ref{33}) we have
\be r_{ab}= ckb^2 \delta_{ab}. \label{34}
\ee
By (\ref{20}), (\ref{30}), (\ref{215}), (\ref{34case1}) and (\ref{34}), we get
$$r_{00} = r_{11} (y^1)^2 + 2 \bar{r}_{10} y^1 + \bar{r}_{00}= kcb^2  \alpha^2 +k(1-c) \beta^2 ,$$
which is equivalent to
\be r_{ij}= kcb^2  a_{ij} + k(1-c) b_i b_j  . \label{35case1}
\ee
Thus we obtain (\ref{lem4.103case1}).

In this case, (\ref{31}) becomes
\be 2 b^3(\Psi s^2+2 \Psi R b^4 + \Pi b^2 s - R b^2)+2 c b^3 \Psi (b^2-s^2)=\mu s^2+\nu(b^2-s^2),\label{36case1}
\ee
where $\mu=\mu(x)$ and $\nu=\nu(x)$  satisfying
$$ k\mu = G^1_{11}-H, \ \ \ k \nu \bar{\alpha}^2 = \bar{G}^1_{00}  .$$
Differentiating (\ref{36case1}) with respect to $s$ for three times yields
\be 2 b^3(\Psi s^2+2 \Psi R b^4 + \Pi b^2 s - R b^2)_{sss}+2 c b^5 \Psi_{sss}=0.\label{37}
\ee
Then $c$ is a function of $b^2$. Therefore, letting $s=0$ in (\ref{36case1}), it is easy to see that $\nu$ is a function of $b^2$. Then by (\ref{36case1})  $\mu$ is a function of $b^2$ too.

Substituting the expressions of $\Psi$ , $\Pi$ and $R$ into (\ref{36case1}) yields
$$\{b^3 [(1-c)s^2+cb^2]+[ ( \nu-\mu )s^2 - \nu b^2](b^2-s^2) \}\phi_{22}=2 b^5  (\phi_1-s \phi_{12})-[(\nu-\mu )s^2 - \nu b^2](\phi-s\phi_2).$$
Thus we obtain (\ref{lem4.102case1}).

\qed
\section{General solutions of (\ref{02})}
In this section, we first give the general solutions of (\ref{02}). Then some special solutions can be constructed. Our method is to take a variable substitution such that (\ref{02}) be a simplier equation.\\
{\bf Proof of Theorem \ref{thm1.2}:}
Since $s$ and $b$ are arbitrary numbers with $|s|\leq b <b_0$, put
\be \psi :=(\phi-s\phi_2)e^{\int{\frac{c-1}{2 b^2}}db^2} \sqrt{b^2-s^2}. \label{38}
\ee
Then
\be \psi_1 =(\phi-s\phi_2)_1 e^{\int{\frac{c-1}{2 b^2}}db^2} \sqrt{b^2-s^2} + (\phi-s\phi_2)e^{\int{\frac{c-1}{2 b^2}}db^2} \Big\{\frac{(c-1) \sqrt{b^2-s^2}}{2 b^2 }+\frac{1}{2 \sqrt{b^2-s^2}}\Big\}, \label{39}
\ee
\be \psi_2 =(\phi-s\phi_2)_2 e^{\int{\frac{c-1}{2 b^2}}db^2} \sqrt{b^2-s^2} - (\phi-s\phi_2) \frac{s }{\sqrt{b^2-s^2}}e^{\int{\frac{c-1}{2 b^2}}db^2}, \label{40}
\ee
where $\psi_1$ is the derivation of $\psi$ with respect to the first variable $b^2$ and $\psi_2$ is the derivation of $\psi$ with respect to the second variable $s$.
It follows from (\ref{39}) and (\ref{40}) that
\be
\label{390}
(\phi-s\phi_2)_1 = \frac{\psi_1 e^{\int{\frac{1-c}{2 b^2}}db^2} }{\sqrt{b^2-s^2}}
-(\phi-s\phi_2) \Big\{\frac{c-1}{2 b^2}+\frac{1}{2 (b^2-s^2)}\Big\},
\ee
\be
\label{400}
(\phi-s\phi_2)_2=  \frac{\psi_2 e^{\int{\frac{1-c}{2 b^2}}db^2} }{\sqrt{b^2-s^2}}
+\frac{s}{b^2-s^2}(\phi-s\phi_2).
\ee
Note that $(\phi-s \phi_2)_2 = - s \phi_{22}$. Substituting (\ref{390}) and (\ref{400}) into (\ref{02}) yields
\be 2 b^5  s \psi_1 + \{b^3 [(1-c)s^2+cb^2]+[ (\nu-\mu)s^2 - \nu b^2](b^2-s^2) \} \psi_2=0. \label{41}
\ee
The characteristic equation of PDE (\ref{41}) is
\be \frac{db^2}{2 b^5  s} = \frac{d s}{b^3 [(1-c)s^2+cb^2]+[ (\nu-\mu)s^2 - \nu b^2](b^2-s^2)},
\label{42}
\ee
which is equivalent to
\be 2s \frac{d s}{db^2} = \frac{(1-c)s^2+cb^2}{ b^2}+ \frac{[(\nu-\mu)s^2- \nu b^2] (b^2-s^2)}{b^5 }.
\label{43}
\ee
Let
\be \chi(b^2) = s^2(b^2) - b^2.
\label{44}
\ee
Differentiating (\ref{44}) with respect to $b^2$ yields
\be \frac{d \chi}{db^2} = 2 s\frac{d s}{db^2} - 1.
\label{45}
\ee
Substituting (\ref{43}) into (\ref{45}) yields
\be
\label{46}
\frac{d \chi}{db^2} = - \frac{(c-1)b-\mu }{b^3} \chi - \frac{\nu-\mu}{b^5} \chi^2,
\ee
which is a Bernoulli equation. It can be rewritten as
$$\frac{d}{db^2}(\frac{1}{\chi}) =   \frac{(c-1)b-\mu }{b^3} \frac{1}{\chi} + \frac{\nu-\mu}{b^5}.$$
By solving this linear  ODE of $\frac{1}{\chi}$, it can be obtained that
\be \frac{1}{\chi} =-e^{  \int{ \frac{(c-1)b-\mu }{b^3} }db^2 } \Big[ \lambda -\int{ \frac{\nu-\mu}{b^5} e^{ \int{ \frac{(1-c)b+\mu }{b^3}}db^2 }}db^2 \Big],
\label{47}
\ee
where $\lambda$ is an arbitrary constant.
Then by (\ref{44}) and (\ref{47}) we get
\be
\label{77}
s=s(b^2)=\sqrt{b^2-\frac{e^{  \int{ \frac{(1-c)b+\mu }{b^3}}db^2 }}{ \lambda -\int{ \frac{\nu-\mu}{b^5} e^{ \int{ \frac{(1-c)b+\mu }{b^3}}db^2 }}db^2}}.
\ee
Obviously, $\chi=0$ is also a solution of (\ref{46}), it is excluded in this situation.
Therefore, it follows from (\ref{77}) that
$$\frac{b^2-s^2}{e^{ \int{ \frac{(1-c)b+\mu }{b^3}}db^2 }  + (b^2-s^2) \int{ \frac{\nu-\mu}{b^5} e^{ \int{ \frac{(1-c)b+\mu }{b^3}}db^2 }}db^2} =  \frac{1}{\lambda}.$$
Thus the solution of (\ref{41}) is
\be
\psi = \Phi \Big(\frac{b^2-s^2}{e^{ \int{ \frac{(1-c)b+\mu }{b^3}}db^2 }  + (b^2-s^2) \int{ \frac{\nu-\mu}{b^5} e^{ \int{ \frac{(1-c)b+\mu }{b^3}}db^2 }}db^2} \Big),
\label{48}
\ee
where $\Phi$ is an arbitrary differentiable function.
Plugging (\ref{48}) into (\ref{38}), we have
\be
\phi-s \phi_2 = \Phi(\zeta(b^2,s)) e^{\int{\frac{1-c}{2 b^2}}db^2} \frac{1}{\sqrt{b^2-s^2}}.
\label{49}
\ee
where $$\zeta(b^2,s)= \frac{b^2-s^2}{e^{ \int{ \frac{(1-c)b+\mu }{b^3}}db^2 }  + (b^2-s^2) \int{ \frac{\nu-\mu}{b^5} e^{ \int{ \frac{(1-c)b+\mu }{b^3}}db^2 }}db^2}, $$
Set
\be
\phi=s \varphi,
\label{50}
\ee
then
\be
\phi-s \phi_2 = - s^2 \varphi_2.
\label{51}
\ee
By (\ref{49}) and (\ref{51}), we get
\be
\varphi= h(b^2)- \xi(b^2) \int{\frac{ \Phi(\zeta(b^2,s))}{s^2 \sqrt{b^2-s^2}}ds}
\label{52}
\ee
for some $C^{\infty}$ functions $h(b^2)$, where $\xi(b^2)= e^{\int{\frac{1-c}{2 b^2}}db^2}.$ \
Hence, plugging (\ref{52}) into (\ref{50}) yields (\ref{04}).

Note that $- s \phi_{22}=(\phi-s \phi_2)_2 $. Therefore, when $\phi$ is given by (\ref{04}),  by Lemma \ref{lem2.1} and (\ref{49}), $F$ is a Finsler metric if and only if  (\ref{lem511}) and (\ref{lem512}) hold.

\qed\\

Then we can construct many Douglas general $(\alpha,\beta)$-metrics $F=\alpha \phi(b^2,\frac{\beta}{\alpha})$ by choosing suitable $c$,  $\mu$, $\nu$ and $\Phi$. The following are some special solutions. Let $\zeta=\zeta(b^2,s)$ is given by (\ref{05}).
\begin{ex}\label{ex1}
Let $\mu=\nu=0$, $c=1$ and $\Phi(\zeta(b^2,s))=\sqrt{\zeta}$, then parts of the solutions of (\ref{04}) are given by
\be
\label{ex11}
\phi(b^2,s)=1+h(b^2) s,
\ee
where $h$ is an arbitrary $C^{\infty}$ function of $b^2$.

In this case, $F=\alpha \phi(b^2,\frac{\beta}{\alpha})$ is a Randers metric, where $\alpha$ and $\beta$ satisfy (\ref{03}).
\end{ex}

\begin{ex}\label{ex2}
 Let $\mu=\nu=0$, $c=1$ and $\Phi(\zeta(b^2,s))=\sqrt{\frac{\zeta}{1-\zeta}}$, then parts of the solutions of (\ref{04}) are given by
\be
\label{ex21}
\phi(b^2,s)=h(b^2) s + \frac{\sqrt{1-b^2+s^2}}{1-b^2},
\ee
where $h$ is an arbitrary $C^{\infty}$ function of $b^2$.

In this case, $F=\alpha \phi(b^2,\frac{\beta}{\alpha})$ is also a Randers metric, where $\alpha$ and $\beta$ satisfy (\ref{03}).
\end{ex}

The metrics in the following example are just the metrics in \cite{LS1}.
\begin{ex}\label{ex3}
Let $\mu=\nu=0$, then the solutions of (\ref{04}) are given by
\be
\label{ex31}
\phi(b^2,s)=h(b^2)s -e^{\int{\frac{1-c}{2 b^2}}db^2} s \int{\frac{ \Phi(\zeta(b^2,s))}{s^2 \sqrt{b^2-s^2}}ds}
\ee
where $\zeta(b^2,s)= (b^2-s^2)e^{\int{\frac{c-1}{b^2}}db^2}$, $h=h(b^2)$ and $c=c(b^2)$ are two scalar functions.

In this case, $F=\alpha \phi(b^2,\frac{\beta}{\alpha})$ is a Douglas metric, where $\alpha$ and $\beta$ satisfy (\ref{03}).
Particularly,
\be
\label{ex32}
\phi(b^2,s)=1+b^{2 c}+h(b^2) s+b^{2(c-1)}s^2
\ee
when $c=const.$ and  $\Phi(\zeta(b^2,s))=(1+\zeta)\sqrt{\zeta}$. In this case, the corresponding general $(\alpha,\beta)$-metrics
are of Douglas type. They are just the special solutions (iii) in \cite{LS1}.
\end{ex}

By choosing different $\mu$ and $\nu$ from the above example, some new Douglas metrics can be found in the following example.
\begin{ex}\label{ex4}
 Let $c=1-b^2>0$, $\mu=\frac{b^5}{1-b^2}$,  and $\nu=\frac{2b^5}{1-b^2}$, then the solutions of (\ref{04}) are given by
\be
\label{ex41}
\phi(b^2,s)=h(b^2)s -\sqrt{e^{b^2}} s \int{\frac{ \Phi(\zeta(b^2,s))}{s^2 \sqrt{b^2-s^2}}ds},
\ee
where $\zeta(b^2,s)= \frac{(b^2-s^2) (1-b^2)}{1+b^2-s^2}$, $h=h(b^2)$ is a arbitrary $C^{\infty}$ function of $b^2$.
By setting $\Phi(\zeta(b^2,s))=\frac{\sqrt{\zeta}}{(1-\zeta)^{3/2}}$ in (\ref{ex41}),
we have
\be
\label{ex42}
\phi(b^2,s)=h(b^2) s + \frac{(1+b^2)(1+b^4-b^2s^2)+s^2(1-b^2)}{(1+b^4)^2} \sqrt{\frac{(1-b^2)e^{b^2}}{1+b^4-b^2s^2}}.
\ee
In this case, $F=\alpha \phi(b^2,\frac{\beta}{\alpha})$ is a new Douglas metric, where $\alpha$ and $\beta$ satisfy (\ref{03}).
\end{ex}

By choosing suitable functions $c(b^2)$, $\mu(b^2)$, $\nu(b^2)$ and $\Phi(\zeta(b^2,s))$, one can construct more new Douglas metrics.

\section{Special solutions of (\ref{03}) }
In this section, we discuss the solutions  of (\ref{03}) in Theorem \ref{thm1.1} in two cases.
Firstly, we give a special solution of (\ref{03}) when $\alpha$ is a projectively flat Riemannian metric. Secondly, we give the special solutions when $\alpha$ is conformal to $|y|$.
In \cite{LS1}, we proved the following fact
\begin{lem}\label{lem7.1}
Let $\widetilde{\beta} = \rho(b^2) \beta$, then
$$\widetilde{b}_{i|j} = \rho b_{i|j} + 2 \rho' b_i (r_j+s_j).$$
\end{lem}

By Beltrami's theorem, projectively flat Riemannian metrics are of constant sectional curvature. Let $\alpha=\sqrt{a_{ij}y^iy^j}$ be a projectively flat Riemannian metric of constant sectional curvature $\kappa$. Then there is a local coordinate system such that
$$\alpha=\frac{\sqrt{(1+\kappa|x|^2)|y|^2 - \kappa \langle x,y\rangle^2 }}{1+\kappa|x|^2}.$$
Then we can give a special solution of (\ref{03}).
\begin{ex}\label{ex7.1}
Let
$$\alpha=\frac{\sqrt{(1+\kappa|x|^2)|y|^2 - \kappa \langle x,y\rangle^2 }}{1+\kappa|x|^2}$$
and
\be
\beta =\frac{b\{\delta_1 \langle x,y\rangle + (1+\kappa|x|^2)\langle a,y\rangle - \kappa \langle a,x\rangle \langle x,y\rangle\} }{\sqrt{(b^2+\delta_2)(1+\kappa|x|^2)^3}},
\label{prop1}
\ee
where $\kappa$ is the  constant sectional curvature of $\alpha$, $\delta_1$ and $\delta_2$ are constant numbers and $a$ is a constant vector. Suppose that $\widetilde{\beta} = \rho(b^2) \beta$ satisfies
\be
\label{prop2}
\widetilde{b}_{i|j}=kcb^2\rho a_{ij}.
\ee
Thus $\widetilde{\beta}$ is a conformal 1-form of $\alpha$. By the result in \cite{ShX}, we get
\be
\label{prop0}
\widetilde{\beta}=\frac{\delta_1 \langle x,y\rangle + (1+\kappa|x|^2)\langle a,y\rangle - \kappa \langle a,x\rangle \langle x,y\rangle }{(1+\kappa|x|^2)^{\frac{3}{2}}},
\ee
\be
\label{prop3}
\widetilde{b}_{i|j} = \frac{\delta_1 - \kappa \langle a,x\rangle}{\sqrt{1+\kappa|x|^2}} a_{ij}.
\ee
Therefore, by (\ref{prop1}) and (\ref{prop0}), we have
\be
\label{prop4}
\rho=\frac{\sqrt{b^2+\delta_2}}{b}.
\ee
Plugging (\ref{prop4}) into (\ref{prop2}) yields
\be
\label{prop5}
\widetilde{b}_{i|j}=kcb\sqrt{b^2+\delta_2} a_{ij}.
\ee
Comparing (\ref{prop3}) with (\ref{prop5}),  we get
\be
\label{prop6}
kc=\frac{\delta_1 - \kappa \langle a,x\rangle}{b\sqrt{(b^2+\delta_2)(1+\kappa|x|^2)}}.
\ee
In this case, by (\ref{prop2}), (\ref{prop4}),(\ref{prop6}) and  Lemma \ref{lem7.1}, we obtain
\begin{equation}
\label{prop8}
\begin{split}
b_{i|j}&= \frac{1}{\rho} \widetilde{b}_{i|j} + 2 \Big(\frac{1}{\rho}\Big)' \widetilde{b_i} \widetilde{r_j}=\frac{1}{\rho}  (kcb^2\rho  a_{ij}) + 2 \Big(\frac{-\rho'}{\rho^2}\Big) \widetilde{b_i} (kcb^2 \rho  \widetilde{b_j}) =kcb^2 a_{ij}-2 \rho' \rho kcb^2 b_i b_j \\
&=\frac{b(\delta_1 - \kappa \langle a,x\rangle)}{\sqrt{(b^2+\delta_2)(1+\kappa|x|^2)}} a_{ij}+\frac{\delta_2(\delta_1 - \kappa \langle a,x\rangle)}{b^3 \sqrt{(b^2+\delta_2)(1+\kappa|x|^2)}}b_i b_j .
\end{split}
\end{equation}
Thus we get a special solution of (\ref{03}) when
$$c=\frac{b^2}{\delta_2+b^2}, \ \
k=\frac{(\delta_2+b^2)(\delta_1 - \kappa \langle a,x\rangle)}{b^3 \sqrt{(b^2+\delta_2)(1+\kappa|x|^2)}}.$$
\end{ex}

To consider the case when $\alpha$ is not projectively flat,  we give two special solutions when $\alpha$ is conformal to $|y|$ in the following example.
\begin{ex}\label{ex7.2}
Let
$$\alpha = \frac{1}{2|x|}|y|, \ \ \beta = 2\varepsilon e^{-|x|^2} \langle x,y\rangle,$$
where $\varepsilon\neq0$ is a constant.
The Christoffel symbols of $\alpha$ are given by
$$\Gamma^i_{jk} = - \frac{1}{|x|^2} \{x_k \delta^i_j + x_j \delta^i_k - x^i \delta_{jk}\}.$$
Then $$b_{i|j} = 4\varepsilon\frac{ 1-  |x|^2}{ |x|^2}e^{-|x|^2} x_i x_j.$$
Thus we get a special solution when $c=0$ and $k =\frac{ 1-  |x|^2}{\varepsilon |x|^2}e^{|x|^2} $.
\end{ex}
\begin{ex}\label{ex7.3}
Let
$$\alpha = \frac{1}{1 + |x|^2}|y|, \ \ \beta = \frac{1 + |x|^2}{1 - |x|^2} \langle x,y\rangle.$$
The Christoffel symbols of $\alpha$ are given by
$$\Gamma^i_{jk} = - \frac{2}{1+|x|^2} \{x_k \delta^i_j + x_j \delta^i_k - x^i \delta_{jk}\}.$$
Then $$b_{i|j} = \delta_{ij} + \frac{4( 2-  |x|^2)}{(1 - |x|^2)^2} x_i x_j.$$
Thus we get a special solution by choosing $kcb^2=(1 + |x|^2)^2$ and $k(1-c) =\frac{4( 2-  |x|^2)}{(1 + |x|^2)^2} $.
\end{ex}

{}


\begin{thebibliography}{1}


\bibitem{B1} S. B\'{a}cs\'{o} and M. Matsumoto, {\it On Finsler spaces of Douglas type. A generalization of the notion of Berwald space}, Publ. Math. Debrecen  51 (1997), 385-406.

\bibitem{C1} B. Chen, Z. Shen and L.Zhao, {\it On a class of Ricci-flat Finsler metrics in Finsler geometry}, J. Geom. Phys. 70 (2013), 30-38.

\bibitem{D1} J. Douglas, {\it The general geometry of paths}, Ann. of Math. 29(1927-1928), 143-168.

\bibitem{L1} B. Li, Y. Shen and Z. Shen, {\it On a class of Douglas metrics}, Stud. Sci. Math. Hung. 46(3) (2009), 355-365.

\bibitem{LS1} B. Li and Z. Shen, {\it On a class of projectively flat Finsler metrics}, Int. J. Math. 27(4) (2016), 1650052.

\bibitem{MST} X. Mo, N. M. Sol\'{o}rzano and K. Tenenblat, {\it On spherically symmetric Finsler metrics with vanishing Douglas curvature}, Diff. Geom. Appl. 31 (2013), 746-758.

\bibitem{R1} G. Randers, {\it On an asymmetric metric in the four-space of general relativity}, Phys. Rev. 59 (1941), 195-199.

\bibitem{S1} E. Sevim, Z. Shen and L. Zhao, {\it On a class of Ricci-flat Douglas metrics}, Int. J. Math. 23(6) (2012), 1250046.

\bibitem{Sh2} Z. Shen, {\it Differential Geometry of Spray and Finsler Spaces}, Kluwer Academic Publishers, 2001.

\bibitem{Sh1} Z. Shen, {\it On Projectively flat $(\alpha,\beta)$-metrics}, Can. Math. Bull. 52(1) (2009), 132-144.

\bibitem{ShX} Z. Shen and H. Xing, {\it On Randers Metrics with Isotropic S-Curvature}, Acta Math. Sinica 24(5) (2008), 789-796.

\bibitem{Y1} C. Yu and H. Zhu, {\it On a new class of Finsler metrics}, Diff. Geo. Appl. 29 (2011), 244-554.

\bibitem{Zh1} H. Zhu, {\it On general $(\alpha,\beta)$-metrics with vanishing Douglas curvature}, Int. J. Math. 26(9) (2015), 1550076.





\end{thebibliography}
\end{document}